\documentclass[leqno]{article}
\usepackage{amsmath}
\usepackage{amssymb,amscd}
\usepackage{verbatim}

\renewcommand{\baselinestretch}{1.2}
\begin{document}
\hyphenation{semi-per-fect Grothen-dieck}
\newtheorem{Lemma}{Lemma}[section]
\newtheorem{Th}[Lemma]{Theorem}
\newtheorem{Prop}[Lemma]{Proposition}
\newtheorem{OP}[Lemma]{Open Problem}
\newtheorem{Cor}[Lemma]{Corollary}
\newtheorem{Fact}[Lemma]{Fact}
\newtheorem{Def}[Lemma]{Definition}
\newtheorem{Not}[Lemma]{Notation}
\newtheorem{Ex}[Lemma]{Example}
\newtheorem{Exs}[Lemma]{Examples}
\newtheorem{Rem}[Lemma]{Remark}
\newenvironment{Remarks}{\noindent {\bf Remarks.}\ }{}
\newtheorem{Remark}[Lemma]{Remark}
\newenvironment{Proof}{\noindent{\sc Proof.}\ }{~\rule{1ex}{1ex}\vspace{0.5truecm}}

\newenvironment{Proofpbinequalities}{\noindent{\sc{ Proof of  Theorem~\ref{pbinequalities}.}}\ }{~\rule{1ex}{1ex}\vspace{0.5truecm}}
\newcommand{\End}{\mbox{\rm End}}
\newcommand{\K}{\mbox{\rm K.dim}}
\newcommand{\Hom}{\mbox{\rm Hom}}
\newcommand{\infsupp}{\mathrm{inf\mbox{-}supp}\,}
\newcommand{\finsupp}{\mathrm{fin\mbox{-}supp}\,}
\newcommand{\Ext}{\mbox{\rm Ext}}
\newcommand{\supp}{\mbox{\rm supp}\,}
\newcommand{\Supp}{\mbox{\rm Supp}\,}
\newcommand{\Max}{\mbox{\rm Max}}
\newcommand{\cl}{\mbox{\rm cl}}
\newcommand{\add}{\mbox{\rm add}}
\newcommand{\Inv}{\mbox{\rm Inv}}
\newcommand{\rk}{\mbox{\rm rk}}
\newcommand{\Tr}{\mbox{\rm Tr}}
\newcommand{\Sat}{\mbox{\bf Sat}}
\newcommand{\card}{\mbox{\rm card}}
\newcommand{\Ann}{\mbox{\rm Ann}}
\newcommand{\proj}{\mbox{\rm proj-}}
\newcommand{\codim}{\mbox{\rm codim}}
\newcommand{\B}{\mathcal{B}}
\newcommand{\Scal}{\mathcal{S}}
\newcommand{\Cong}{\mbox{\rm Cong}}
\newcommand{\Spec}{\mbox{\rm Spec-}}
\newcommand{\coker}{\mbox{\rm coker}}
\newcommand{\Cl}{\mbox{\rm Cl}}
\newcommand{\Ses}{\mbox{\rm Ses}}
\newcommand{\im}{\mbox{\rm Im}}
\newcommand{\Cal}[1]{{\cal #1}}
\newcommand{\+}{\oplus}
\newcommand{\N}{\mathbb N}
\newcommand{\No}{{\mathbb N}_0}
\newcommand{\Z}{\mathbb{Z}}
\newcommand{\Q}{\mathbb{Q}}
\newcommand{\C}{\mathbb{C}}
\newcommand{\T}{\mathbb{T}}
\newcommand{\R}{\mathbb{R}}
\newcommand{\notsim}{\sim\makebox[0\width][r]{$\slash\;\,$}}
\newcommand{\Mod}{\mbox{\rm Mod-}}
\newcommand{\lmod}{\mbox{\rm -mod}}
\newcommand{\mspec}{\mbox{\rm Max}}
\renewcommand{\dim}{\mathrm{dim}\,}

\title {Infinitely generated projective modules over pullbacks of rings}

\author{Dolors Herbera\thanks{The final version of this paper was written while the author was visiting NTNU (Tondheim, Norway), she thats her host for the kind hospitality.
She was also partially supported by MEC-DGESIC (Spain) through
Project MTM2008--06201-C02-01, and by the Comissionat Per
Universitats i Recerca de la Generalitat de Catalunya through
Project 2005SGR00206.}
\\
 Departament de Matem\`atiques, \\
Universitat Aut\`onoma de Barcelona, \\ 08193 Bellaterra
(Barcelona), Spain\\ e-mail: dolors@mat.uab.cat  \and Pavel P\v
r\'\i hoda\thanks{Supported by GA\v CR 201/09/0816  and research project MSM 0021620839.} \\
Charles University, Faculty of Mathematics and Physics, \\Department
of Algebra, Sokolovsk\'a 83,\\ 18675 Praha 8, Czech Republic\\
e-mail: prihoda@karlin.mff.cuni.cz}

\date{\phantom{ciao}}

\maketitle

\begin{abstract} We use pullbacks of rings to realize the submonoids $M$ of
$(\N _0\cup\{\infty\})^k$ which are the set of solutions of a finite
system of linear diophantine inequalities as the monoid of
isomorphism classes of countably generated projective right
$R$-modules over a suitable semilocal ring. For these rings, the
behavior of countably generated projective left $R$-modules is
determined by the monoid $D(M)$ defined by  reversing the
inequalities determining the monoid $M$. These two monoids are not
isomorphic in general. As a consequence of our results we show that
there are semilocal rings such that all its projective right modules
are free but this fails for   projective left modules. This answers
in  the negative a question posed by Fuller and Shutters \cite{FS}.
We also provide a rich variety of examples of semilocal rings having
non finitely generated projective modules that are finitely
generated modulo the Jacobson radical.
\end{abstract}

After the paper of Bass \cite{bass}  there seemed to be the general
belief that the theory of infinitely generated projective modules
\emph{invited little interest}. However some  of the developments in
the representation theory of finite dimensional algebras
\cite{ringel} and subsequent ones in integral representation theory
have drawn the attention to the infinite dimensional representations
\cite{Rump}, \cite{BCK}. Also the study of the direct sum
decomposition of infinite direct sums of modules over general rings
requires a good knowledge of the behavior of all projective modules
\cite{Pun1}. As a result of this pressure, interesting general
theory on projective modules has recently appeared \cite{P1},
\cite{P2} and it has been shown that examples of rings such that not
all projective modules are direct sum of finitely generated are
relatively frequent \cite{PP} and the behavior can be quite complex
even for noetherian rings \cite{FH}. In this paper we continue this
line of work by providing further examples of such rings. All of
them are semilocal rings, that is, rings that are semisimple
artinian modulo the Jacobson radical.

Our study makes essential use of the result proved by P. P\v r\'\i
hoda in \cite{P1} that, over an arbitrary ring, projective modules
are isomorphic if and only if they are isomorphic modulo the
Jacobson radical. For a semilocal ring $R$ this implies that the
monoid of isomorphism classes of countably generated projective
right (or left) $R$-modules can be seen as a submonoid of $(\No\cup
\{\infty\})^k$ for a suitable $k\ge 1$, cf. \S \ref{preliminaries}
for the precise definitions.

In \cite{HP}, we characterized the class of monoids that can be
realized as monoid of isomorphism classes of countably generated
projective right (or left) modules over a noetherian semilocal ring
as \emph{essentially} the set of solutions in $\No\cup \{\infty\}$
of finite homogeneous systems of diophantine linear equations. In
Theorem \ref{pbinequalities} we show that any monoid $M$ which is
the set of solutions in $\No\cup \{\infty\}$ of a finite homogeneous
system of diophantine linear inequalities can   also be realized as
monoid of isomorphism classes of countably generated projective
right modules over a suitable semilocal ring $R$. In the examples we
construct, the monoid of isomorphism classes of countably generated
projective left $R$-modules  is the set of solutions in $\No\cup
\{\infty\}$ of the system obtained by reversing the inequalities of
the system defining $M$. While in the noetherian case the monoid of
countably generated  projective right modules is isomorphic to the
one of countably generated projective left modules, as we show in
this paper, this is no longer true for general semilocal rings.

In this paper we emphasize in the study of projective modules that
are not finitely generated but that they are finitely generated
modulo the Jacobson radical. The first example of this kind was
provided by Gerasimov and Sakhaev in \cite{GS}, and the construction
was further developed by Sakhaev in \cite{Sa}. Other examples appear
when studying the direct sum decomposition of infinite direct sums
of uniserial modules \cite{Pun1}, \cite{FHS} and \cite{P1}. From
these examples it seemed that the existence of such projective
modules is rare and very difficult to handle. With our methods we
can produce a wide variety of examples where such projectives exist
and where their behavior is \emph{under control}. In our examples,
the countably generated projective modules that are finitely
generated modulo the Jacobson radical, correspond to the solutions
in $\No$ of the system of inequalities. Between them we distinguish
the finitely generated ones as the ones that fulfill the equality.

The techniques we use in this paper are an extension of the ones in
\cite{HP}. As the title indicates, our rings are constructed as
pullbacks of suitable rings, and we take advantage of
\cite[Theorems~2.1, 2.2 and 2.3]{milnor} in which Milnor describes
\emph{all} projective modules over a class of  ring pullbacks. A key
ingredient  is the Gerasimov-Sakhaev example mentioned above and the
computation of its monoid of isomorphism classes of countably
generated projective right (and left) modules done in \cite{DPP}.

In \S \ref{preliminaries} we  give an overview of the paper: we
introduce the monoids of projective modules, we define in a precise
way the class of monoids that we will realize in section
\ref{realizing} as monoids of  countably generated projective right
modules and of  countably generated projective left modules over
suitable semilocal rings, and we state our main realization Theorem
\ref{pbinequalities}.

In section \ref{pure} we develop some   theory on projective modules
that are finitely generated modulo the Jacobson radical which
essentially follows \cite{Sa}. Theorem \ref{flatproj} is a slight
generalization of the main result in \cite{plans}.

In section \ref{examples} we compute some particular examples to
illustrate the consequences of Theorem~\ref{pbinequalities}. For
instance, in \ref{noniso}, we construct a semilocal ring such that
all projective left $R$-modules are free while $R$ has a nonzero
(infinitely generated) right projective module that is not a
generator. Such an example also shows that the notion of p-connected
ring is not left-right symmetric; this answers in the negative a
question in \cite[page 310]{FS}. Recall that, following Bass
\cite{bass}, a ring is (left) $p$-connected if every nonzero left
projective module is a generator.

We also provide examples showing that if $R$ is a semilocal ring
such that $R/J(R)\cong D_1\times D_2$ and $R$ has a countably
generated, but not finitely generated, projective module that is
finitely generated modulo the Jacobson radical then there is still
room for countably generated (right and left, or just right)
projective modules that are not direct sums of projective modules
that are finitely generated modulo the Jacobson radical. This answers
in the negative a question formulated in \cite[page 3261]{DPP}.

In section \ref{monoids} we develop some properties of the monoids
defined by inequalities. Finally, in section \ref{realizing} we
prove Theorem \ref{pbinequalities}.

\section{Preliminaries and overview} \label{preliminaries}

All our rings are associative with $1$, and ring morphism means
unital ring morphism.

\subsection{Monoids of projective modules}
Let $R$ be a ring. Let $V^*(R_R)=V^*(R)$ ($V^*({}_RR)$) be the set
of isomorphism classes of countably generated projective right
(left) $R$-modules.  If $P$ and $Q$ are countably generated
projective right $R$-modules then the direct sum induces an addition
on $V^*(R)$ by setting $\langle P\rangle +\langle Q\rangle =\langle
P\oplus Q\rangle$, so that $V^* (R)$ is an additive monoid.
Similarly, $V^*({}_RR)$ is also an additive monoid.

Let $V(R)$ be the set of isomorphism classes of finitely generated
right (or left) $R$-modules. Again $V(R)$ is an additive monoid,
which can be identified with a submonoid of $V^*(R)$. Since the
functor $\mathrm{Hom}_R(-,R)$ induces a duality between the category
of finitely generated projective right $R$-modules and the category
of finitely generated projective left $R$-modules we identify
$V({}_RR)$ with $V(R)$. So that, we also see $V(R)$ as a submonoid
of $V^*({}_RR)$.

Another interesting submonoid of $V^*(R)$ is $W(R_R)=W(R)$ which we
define as the set of isomorphism classes of countably generated
projective right $R$-modules that are pure submodules of finitely
generated projective modules. The submonoid of $V^*({}_RR)$,
$W({}_RR)$ is defined in a similar way. Clearly, $V(R)\subseteq
W(R)\subseteq V^*(R)$, and $V(R)\subseteq W({}_RR)\subseteq
V^*({}_RR)$. Notice that $W(R)\setminus V(R)$ is also a semigroup.

Along the paper we will find many examples of (semilocal) rings $R$
with non trivial $W(R)$. Now we give a different kind of example.

\begin{Ex} \emph{\cite{bass}} Let $R$ denote the ring of continuous real valued
functions over the interval $[0,1]$. Let
\[I=\{f\in R\mid \mbox{ there exists $\varepsilon >0$ such that }f([0,\varepsilon])=0\}\]
then $I$ is a projective pure ideal of $R$, cf. \cite[Example
3.3]{FHS} or \cite[p. 3263]{DPP}.\end{Ex}

The notation   $W(R)$ is borrowed from the $C^*$-algebra world, as
we think on this monoid as an algebraic analogue of the Cuntz monoid
defined in $C^*$-algebras.

\subsection{The semilocal case}\label{semilocal}  A ring $R$ is
said to be semilocal if modulo its Jacobson radical $J(R)$ is
semisimple artinian, that is, $R/J(R)\cong M_{n_1}(D_1)\times \dots
\times M_{n_k}(D_k)$ for suitable division rings $D_1,\dots ,D_k$.
For the rest of our discussion we  fix an onto ring homomorphism
$\varphi \colon R\to  M_{n_1}(D_1)\times \dots \times M_{n_k}(D_k)$
such that $\mathrm{Ker}\, \varphi =J(R)$.

Let $V_1,\dots ,V_k$ denote a fixed ordered set of representatives
of the isomorphism classes of simple right $R$-modules such that
$\mathrm{End}_R(V_i)\cong D_i$. Let us also fix $W_1,\dots ,W_k$,
where $W_i=\mathrm{Hom}_R(V_i,R/J(R))$ for $i=1,\dots ,k$, as an
ordered set of representatives of simple left $R$-modules.

If $P_R$ is a countably generated projective right $R$-module then
$P/PJ(R)\cong V_1^{(I_1)}\oplus \cdots \oplus V_k^{(I_k)}$ and the
cardinality of the sets $I_1,\dots ,I_k$ determines the isomorphism
class of $P/PJ(R)$. By \cite{P1} (cf. Theorem \ref{Pavel})
projective modules are determined, up to isomorphism, by its
quotient modulo the Jacobson radical. So that, for a semilocal ring
$R$, to describe $V^*(R)$ we only need to record the cardinality of
the sets $I_i$ for $i=1,\dots,k$. A similar situation holds for
projective left $R$-modules.

Note that, by Theorem \ref{Pavel}(i), in the case of semilocal rings
$$W(R)=\{\langle P\rangle\in V^*(R) \mid P/PJ(R)\mbox{ is finitely
generated}\}.$$ Similarly, for $W({}_RR)$.

\subsection{The dimension monoids for semilocal rings} Let $\N=\{1,2,\dots\}$ and $\No =\N \cup \{0\}$. We also
consider the monoid $\No^*=\No \cup \{\infty\}$ with the addition
determined by the addition on $\No$ extended by the rule
$n+\infty=\infty +n=\infty$ for any $n\in \No^*$.

Following the notation of \S \ref{semilocal}, if $P$ is a countably
generated projective right $R$-module such that $P/PJ(R)\cong
V_1^{(I_1)}\oplus \cdots \oplus V_k^{(I_k)}$ we set $\dim _\varphi
(\langle P\rangle))=(m_1,\dots ,m_k)\in (\No^*)^k$ where, for
$i=1,\dots,k$, $m_i=\vert I_i\vert$ if $I_i$ is finite and
$m_i=\infty$ if $I_i$ is infinite. Therefore $\dim _\varphi\colon
V^*(R)\to (\No^*)^k$ is a monoid morphism. Similarly, we define a
monoid morphism $\dim _\varphi\colon V^*({}_RR)\to (\No^*)^k$.

By Theorem \ref{Pavel}(ii), $\mathrm{dim}_\varphi \colon
V^*(R)\to(\No^*)^k$ and $\mathrm{dim}_\varphi \colon
V^*({}_RR)\to(\No^*)^k$ are monoid monomorphisms. Note that
$\dim_\varphi (\langle R\rangle)=(n_1,\dots ,n_k)\in \N ^k$ and that
$\dim_\varphi (W(R))=\No^k\cap \dim_\varphi (V^*(R))$ while
$\dim_\varphi (W({}_RR))=\No^k\cap \dim_\varphi (V^*({}_RR)).$

\begin{Def} \label{fullaffine} A submonoid $A$ of $\No ^k$ is said
to be full affine if whenever $a$, $b\in A$ are such that $a=b+c$
for some $c\in\No ^k$ then $c\in A$.\end{Def}

The class of full affine submonoids of $\No ^k$ containing an
element $(n_1,\dots ,n_k)\in\N ^k$ is the precise class of monoids
that can be realized as $\dim _\varphi (V(R))$ for a semilocal ring
$R$ such that $\dim _\varphi(\langle R\rangle )=(n_1,\dots ,n_k)$
\cite{FH}.

The general problem we are  interested in is determining which
submonoids of $(\No^*)^k$ can be realized as dimension monoids, that
is, as $\mathrm{dim}_\varphi (V^*(R))$ for a suitable semilocal ring
$R$. We do not know the complete solution of this problem but in the
next definition we single out some classes of monoids  that can be
realized as dimension monoids of  semilocal ring.

\begin{Def} \label{defequations} Let $k\ge 1$.
\begin{itemize}
\item[(i)] A submonoid $M$ of $(\No^*)^k$ is said to be a monoid
\emph{defined by a system of equations} if it  is the set of
solutions in $(\No^*)^k$ of a system of the form
\[D\left(\begin{array}{c}t_1\\\vdots \\ t_k\end{array}\right)\in \left(\begin{array}{c}m_1\No^*\\
\vdots \\ m_n\No^*\end{array} \right) \qquad(*)  \qquad \mbox{and}
\qquad
E_1\left(\begin{array}{c}t_1\\\vdots \\
t_k\end{array}\right)= E_2\left(\begin{array}{c}t_1\\\vdots \\
t_k\end{array}\right)\qquad(**) \] where $D\in M_{n\times k}(\No)$,
$E_1,$ $E_2\in M_{\ell \times k}(\No)$, $m_1,\dots ,m_n\in \N$,
$m_i\ge 2$ for any $i\in \{1,\dots ,n\}$ and $\ell$, $n\ge 0$.

\item[(ii)]  A submonoid $M$ of $(\No^*)^k$ is said to be a monoid \emph{defined
by a system of inequalities} provided that there exist $D\in
M_{n\times k}(\No)$, $E_1,$ $E_2\in M_{\ell \times k}(\No)$, $\ell$,
$n\ge 0$, and $m_1,\dots ,m_n\in \N$ , $m_i\ge 2$ for any $i\in
\{1,\dots ,n\}$, such that $M$ is the set of   solutions in
$(\No^*)^k$ of
\[D\left(\begin{array}{c}t_1\\\vdots \\ t_k\end{array}\right)\in \left(\begin{array}{c}m_1\No^*\\
\vdots \\ m_n\No^*\end{array} \right) \qquad \mbox{and} \qquad
E_1\left(\begin{array}{c}t_1\\\vdots \\
t_k\end{array}\right)\ge E_2\left(\begin{array}{c}t_1\\\vdots \\
t_k\end{array}\right). \]

\item[(iii)] If $M\le (\No^*)^k$ is defined by a system of inequalities as in $(ii)$ we define its dual monoid $D(M)$
as  the set of solutions in $(\No^*)^k$ of
\[D\left(\begin{array}{c}t_1\\\vdots \\ t_k\end{array}\right)\in \left(\begin{array}{c}m_1\No^*\\
\vdots \\ m_n\No^*\end{array} \right)  \qquad \mbox{and} \qquad
E_1\left(\begin{array}{c}t_1\\\vdots \\
t_k\end{array}\right)\le E_2\left(\begin{array}{c}t_1\\\vdots \\
t_k\end{array}\right) \] \end{itemize}
\end{Def}

\begin{Rem} \label{allfg} 1) It is important to notice that $\No ^*$ is no longer a
cancellative monoid. So that, for example, the set of solutions in
$(\No^*)^2$ of the equation $x=y$ is not the same as the set of
solutions of $2x=y+x$.

2) If $M$ is a monoid defined by a system of inequalities then the
monoid $D(M)$ depends on the particular system fixed to define $M$.
For an easy example see Examples \ref{noniso}(ii) and (iii).

3) Let  $A$ be a  submonoid of $\No^k$ containing $(n_1,\dots,
n_k)\in \N ^k$. It was observed by Hochster that $A$ is full affine
if and only if $A$ is the set of solutions in $\No^k$ of a system of
the type appearing in Definition \ref{defequations}(i)(cf. \cite[\S
6]{HP}).

In this case, the monoid $M=A+\infty\cdot A$ is a submonoid of
$(\No^*)^k$ defined by a system of equations \cite[Corollary
7.9]{HP}.
\end{Rem}

\subsection{Realization results. Main result} For further quoting we recall the main result in \cite{HP}
which characterized the monoids $M$ that can be realized as $V^*(R)$
for a semilocal noetherian ring $R$. For this class of rings a
projective module that is finitely generated modulo $J(R)$ must be
finitely generated so that $W(R)=V(R)$ (see, for example,
Proposition \ref{sequence}), and also, by \cite{P2},
$V^*({}_RR)\cong V^*(R)$.

\begin{Th} \label{main} Let $k\in \N$. Let $M$ be a submonoid of $(\No^*)^k$
containing  $(n_1,\dots ,n_k)\in \N ^k$. Then the following
statements are equivalent:
\begin{itemize}
\item[(1)] $M$ is a monoid defined by a system of equations.
\item[(2)] There exist a noetherian semilocal
ring $R$, a semisimple ring $S=M_{n_1}(D_1)\times \dots \times
M_{n_k}(D_k)$, where $D_1,\dots,D_k$ are division rings, and an onto
ring morphism $\varphi \colon R\to S$ with $\mathrm{Ker}\, \varphi
=J(R)$ such that $\dim _\varphi  V^*(R) =M$. Therefore, $\dim
_\varphi  V(R) =M\cap \No^k$.
\end{itemize}
In the above statement, if $F$ denotes a field, $R$ can be
constructed to be an $F$-algebra such that $D_1=\cdots =D_k=E$ is a
field extension of $F$.
\end{Th}

In this paper we shall prove the following realization result

\begin{Th}\label{pbinequalities} Let $k\ge 1$, and let $F$ be a field. Let $M$ be a submonoid of  $(\N _0^*)^k$
defined by a system of  inequalities. Let $D(M)$ denote its dual
monoid. Assume that $M\cap D(M)$ contains an element $(n_1,\dots
,n_k)\in \N ^k$. Then there exist a semilocal $F$-algebra $R$, a
semisimple $F$-algebra $S=M_{n_1}(E)\times \dots \times M_{n_k}(E)$,
where $E$ is a suitable field extension of $F$,  and an onto
morphism of $F$-algebras $\varphi \colon R\to S$ with
$\mathrm{Ker}\, \varphi =J(R)$ satisfying that $\dim _\varphi
V^*(R_R)=M$ and $\dim _\varphi V^*({}_RR)=D(M)$.

Moreover, $\dim _\varphi W(R_R)=M\cap \No^k$, $\dim _\varphi
W({}_RR)=D(M)\cap \No^k$, and $\dim _\varphi V(R)=M\cap D(M)\cap
\No^k$.
\end{Th}

For any semilocal ring $V(R)$ is a finitely generated monoid, so is
$V^*(R)$ for $R$ noetherian and semilocal. As we will show in \S
\ref{monoids}, monoids defined by a system of inequalities are still
finitely generated. But, in general, we do not know whether a monoid
that can be realized as $V^*(R)$ for some semilocal ring $R$ must be
finitely generated.

\section{Projective modules, monoids of projectives and Jacobson radical}\label{pure}

In this section we  want to explain the relation between $W(R_R)$
and $W({}_RR)$ completing the results in \cite{plans}. We also take
the opportunity to state in a (too) precise way results on lifting
maps between projective modules modulo an ideal contained in the
Jacobson radical.

Let $I$ be a two-sided ideal of a ring $R$, let $M$ and $N$ be right
$R$-modules, and let $f\colon M\to N$ denote a module homomorphism.
By the induced homomorphism $\overline{f}\colon M/MI\to N/NI$ we
mean the map defined by $\overline{f}(m+MI)=f(m)+NI$ for any $m\in
M$.

Recall the following well known result.

\begin{Lemma}\label{isofg} Let $R$ be any
ring,   and let $I\subseteq J(R)$ be a two-sided ideal of $R$. Let
$f\colon P\to Q$ be a morphism between finitely generated projective
right $R$-modules. Then $f$ is an isomorphism if and only if the
induced homomorphism $\overline{f}\colon P/PI\to Q/QI$ is an
isomorphism.
\end{Lemma}

In contrast, for general projective modules we have.

\begin{Th}\label{Pavel}  Let $R$ be any
ring,  let $P$ and $Q$ be projective right $R$-modules, and let
$I\subseteq J(R)$ be a two-sided ideal of $R$.
\begin{itemize}
\item[(i)] {\rm \cite[Proposition~6.1]{plans}} A module homomorphism $f\colon P\to Q$
is a pure monomorphism if and only if so is the induced map
$\overline{f}\colon P/PI\to Q/QI$.

\item[(ii)] {\rm \cite[Theorem~2.3 and its proof]{P1}} Let $\alpha \colon P/PI\to Q/QI$ be an isomorphism of right
$R/I$-modules. Let $f\colon P\to Q$ be a module homomorphism such
that $\overline{f}=\alpha$, and let $X$ be a finite subset of $P$.
Then there exists an isomorphism $g\colon P\to Q$ such that
$\overline{g}=\alpha$ and $g(x)=f(x)$ for any $x\in X$.

In particular, $P$ and $Q$ are isomorphic if and only if they are
isomorphic modulo the Jacobson radical.
\end{itemize}
\end{Th}

For further applications we note the following corollary of Theorem
\ref{Pavel}.

\begin{Cor}\label{unit} Let $R$ be a ring, and let $I\subseteq J(R)$ be a two-sided ideal. Let $P$ be a countably generated
projective right $R$-module. Let $f\colon P\to P$ be a homomorphism
such that the induced map $\overline{f}\colon P/PI\to P/PI$ is the
identity, and let $X$ be a finite subset of $P$. Then there exists a
bijective homomorphism $h\colon P\to P$ such that the induced
homomorphism $\overline h=\mathrm{Id}_{P/PI}$ and such that
$hf(x)=x$ for any $x\in X$.
\end{Cor}

\begin{Proof} By Theorem \ref{Pavel}(ii), there exists an
isomorphism $g\colon P\to P$ such that $\overline g
=\mathrm{Id}_{P/PI}$ and $g(x)=f(x)$ for any $x\in X$. Set
$h=g^{-1}$ to conclude.
\end{Proof}

\begin{Lemma}\label{idempotent} Let $R$ be a ring, let $P$ and $Q$
be projective right $R$-modules. Let $I$ be a two-sided ideal of $R$
contained in $J(R)$, and let $\alpha\colon Q/QI\to P/PI$  and
$\beta\colon P/PI\to Q/QI$ be homomorphisms such that $\beta\circ
\alpha =\mathrm{Id}_{Q/QI}$. Let $f\colon Q\to P$ and $g\colon P\to Q$ be
module homomorphisms such that $\overline f=\alpha$ and $\overline
g=\beta$.

If $f\circ g$ is idempotent then $P\cong Q\oplus Q'$ and
$Q'/Q'I\cong (\mathrm{Id}_{P/PI}-\alpha \beta)(P/PI)$.
\end{Lemma}

\begin{Proof} Since $fg(P)$ is a direct summand of $P$,
\[fg(P)/fg(P)I=fg(P)/(fg(P)\cap PI)\cong
\left(fg(P)+PI\right)/PI.\] Since, for any $x\in P$,
$\beta(fg(x)+PI)=\beta (x+PI)$ we deduce that $\beta \colon
fg(P)/fg(P)I\to Q/QI$ is bijective. By Theorem \ref{Pavel}, we
conclude that $Q\cong fg(P)$.

Since $\left( (\mathrm{Id}_P-fg)(P)+PI\right)
/PI=(\mathrm{Id}_{P/PI}-\alpha \beta)\left(P/PI\right)$, it follows
that $Q'=(\mathrm{Id}_P-fg)(P)$ has the claimed properties.
\end{Proof}

\begin{Cor} \label{full} Let $R$ be a ring with Jacobson radical $J(R)$. Let $I\subseteq J(R)$ be a two-sided ideal. Let
$P$ and $Q$ be  projective right $R$-modules such that $Q$ is
finitely generated. If there exists a   projective right
$R/I$-module $X$ such that $P/PI\cong Q/QI\oplus X$ then there
exists a projective right $R$-module $Q'$ such that $P\cong Q\oplus
Q'$ and $Q'/Q'I\cong X$.
\end{Cor}

\begin{Proof} Since $Q$ is finitely generated, the split exact sequence of $R/I$-modules
\[0\to X\to P/PI\stackrel{\beta}{\to} Q/QI\to 0\]
lifts to a split exact sequence
\[0\to \mathrm{Ker}\, g\to P\stackrel{g}{\to} Q\to 0\]
where $\overline{g}=\beta$. Therefore $P\cong Q\oplus \mathrm{Ker}\,
g$. We want to show that $\mathrm{Ker}\, g/(\mathrm{Ker}\, g)I\cong
X$.

Let $\alpha \colon Q/QI\to P/PI$ be such that $\beta \alpha =\mathrm
{Id}_{Q/QI}$, and let $f\colon Q\to P$ be such that $\overline
f=\alpha$. Since $Q$ is finitely generated and $\overline{gf}=\beta
\alpha=\mathrm {Id}_{Q/QI}$, $gf\colon Q\to Q$ is invertible (cf.
Lemma \ref{isofg}). So that, there exists an invertible endomorphism
$h$ of $Q$ satisfying that $\overline{h}=\mathrm{Id}_{Q/QI}$, and
such that $g(fh)=\mathrm{Id}$. Therefore, $(fh)g$ is an idempotent
endomorphism of $P$ and since $(\mathrm{Id}-(fh)g)P=\mathrm{Ker}\,
g$ we conclude, by the second part of Lemma \ref{idempotent}, that
$Q'=\mathrm{Ker}\, g$ has the claimed properties.
\end{Proof}

In the following lemma we recall the properties of sequences
$\{f_n\}_{n\ge 1}$ satisfying that $f_{n+1}f_n=f_n$. Lazard in
\cite{lazard} realized the importance of them to describe pure
ideals of a ring. They play a fundamental r\^ole in constructing
finitely generated flat modules over  semilocal rings that are not
projective or, equivalently, in constructing non-finitely generated
projective modules that are finitely generated modulo the Jacobson
radical.

They were very well analyzed by Sakhaev in several papers, see for
example \cite{Sa}. Recently, they have been extensively re-studied
\cite{FHS}, \cite{plans} and \cite{DPP}.

\begin{Lemma} \label{countablenet} Let $R$ be any ring. Let $P$ be
a right $R$-module and let $f_1,\dots ,f_n, \dots$ be a sequence of
endomorphisms of $P$ satisfying that, for each $n\ge 1$,
$f_{n+1}f_n=f_n$ then,
\begin{itemize}
\item[(i)] $\bigcup _{n\ge 1}f_n\cdot \mathrm{End}_R(P)$ is a projective
pure right ideal of $\mathrm{End}_R(P)$.
\item[(ii)] $Q=\bigcup _{n\ge 1}f_n(P)$ is a pure submodule of $P$  isomorphic to a direct summand of  $P^{(\N)}$. In particular,
if $P$ is projective then so is $Q$.
\end{itemize}
\end{Lemma}

\begin{Proof} $(i).$ This is due to Lazard \cite{lazard}.

$(ii)$. The purity of $I$ inside $S$ gives $I \otimes_{S} P
\hookrightarrow S \otimes_{S} P$. Using the identification $S
\otimes_{S} P \simeq P$, we get $\bigcup _{n\ge 1}f_n(P) \simeq I
\otimes_{S} P$.
Hence the
purity of $Q$ inside $P$ follows from the associativity of the
tensor product and $(i)$.

Consider the countable direct system
\[P_1\stackrel{f_1}{\to}P_2\cdots P_n\stackrel{f_n}{\to}P_{n+1}\cdots\]
where $P=P_n$ for any $n\ge 1$. Since $f_{n+1}f_n=f_n$, the sequence
$\{f_n\}_{n\ge 1}$ induces an injective map $f\colon \varinjlim
P_n\to P$ such that $\mathrm{Im}\, f=Q$. Therefore, $Q$ fits into
the (pure) exact sequence
\[0\to \oplus_{n\ge1} P_n\stackrel{\Phi}{\to}\oplus_{n\ge 1} P_n\to Q\to 0\]
where, for each $n\ge 1$ and letting $\varepsilon _n\colon
P_n\to\oplus_{n\ge 1} P_n$ denote the canonical embedding, the map
$\Phi$ is determined by $\Phi \varepsilon_n (x)= \varepsilon_n
(x)-\varepsilon _{n+1}f_n(x)$ for each $x\in P_n$.

The properties of the sequence of maps $\{f_n\}_{n\ge 1}$ imply that
$\Phi$ splits see, for example, \cite[Proposition
2.1]{angelerisaorin}.
\end{Proof}

\begin{Prop} \label{sequence} Let $R$ be a ring. Let $P_R$ and $Q_R$ be projective
right $R$-modules such that $P_R$ is finite generated. Let $\alpha
\colon Q/QJ(R)\to P/PJ(R)$ and $\beta \colon P/PJ(R)\to Q/QJ(R)$ be
such that $\beta \alpha =\mathrm{Id}_{Q/QJ(R)}$. Let $\varepsilon
\colon Q\to P$ be any module homomorphism such that
$\overline{\varepsilon}=\alpha$. Then there exists a sequence
$f_1,\dots ,f_n,\dots $ of endomorphisms of $P$ such that, for each
$n\ge 1$, $f_{n+1}f_n=f_n$, $\overline{f_n}=\alpha \circ \beta$ and
$Q\cong \varepsilon (Q)=\bigcup _{n\ge 1}f_n(P)$.

Moreover $Q$ is finitely generated if and only if there exists $n_0$
such that $f_{n_0}^2=f_{n_0}$. In this case, $f_{n_0+k}^2=f_{n_0+k}$
for any $k\ge 0$.
\end{Prop}

\begin{Proof}  Let $\varphi \colon P\to
Q$ be a lifting of $\beta$.

Note that $Q_R$ must be a countably generated projective module, so
that we can fix an ascending chain $\emptyset=X_1\subseteq
X_2\subseteq X_3\subseteq \dots \subseteq X_n\subseteq \dots $ of
finite subsets of $Q$ such that $X=\bigcup _{n\ge 1}X_n$ generates
$Q$.

Since $P$ is finitely generated and using Corollary \ref{unit}, we
can construct, inductively, a sequence $\mathrm{Id}_Q=h_1,\dots
,h_n,\dots $ of (auto)morphisms of $Q$ such that if, for each $n\ge
1$, we set $f_n=\varepsilon h_nh_{n-1}\cdots h_1\varphi$ then
$h_{n+1}h_n\cdots h_1\varphi f_n=h_n\cdots h_1\varphi$ and
$h_{n+1}h_n\cdots h_1\varphi \varepsilon (x)=x$ for any $x\in
X_{n+1}$. It can be easily checked  that the homomorphisms
$\{f_n\}_{n\ge 1}$ satisfy the desired properties.

If $Q$ is finitely generated there exists $n_0$ such that
$\varepsilon(Q) = f_{n_0-1}(P)$. Observe that $f_{n_0}f_{n_0-1}=
f_{n_0-1}$ says $f_{n_0}(x) = x$ for any $x \in Q$. In particular,
$f_{n_0+k}^2 = f_{n_0+k}$ for any $k \in \mathbb{N}$.

Conversely, in view of Lemma \ref{idempotent}, if there exists $n_0$
such that $f_{n_0}^2=f_{n_0}$ then $Q$ is isomorphic to $f_{n_0}(P)$
which is a direct summand of $P$. In particular, $Q$ is finitely
generated and $f_{n_0}(P)=f_{n_0+k}(P)$ for any $k\ge 0$. Since
$f_{n_0}$ is idempotent, for any $k\ge 0$,
$f_{n_0+k}=f_{n_0}f_{n_0+k}$ so that
$f_{n_0+k}^2=f_{n_0+k}f_{n_0}f_{n_0+k}=f_{n_0+k}$.
\end{Proof}

\begin{Rem}\label{left} In the situation of Proposition
\ref{sequence}, fix $n\ge 1$. Notice that
$(f_{n+1}-f_n)f_n=f_n-f_n^2$. Since
$\overline{f_{n+1}-f_n}=\overline 0\in \mathrm{End}_R(P/PI)$ and $P$
is a finitely generated projective module,
$u=\mathrm{Id}_P-(f_{n+1}-f_n)$ is a unit such that $uf_n=f_n^2$.

For any $m\in \mathbb{Z}$, set $g_m= u^{-(m+1)}f_nu^m\in
\mathrm{End}_R(P)$. It easily follows that, for any $m\in \Z$,
$g_{m+1}g_m=g_m$ and also that
$(\mathrm{Id}_P-g_{m+1})(\mathrm{Id}_P-g_{m})=\mathrm{Id}_P-g_{m+1}$
so that, by Lemma \ref{countablenet}, $P'_n=\bigcup _{m\ge 0}g_mP$
is a projective pure submodule of $P$ and $Q'_n=\bigcup _{m\le
0}\mathrm{Hom}_R (P,R)(\mathrm{Id}_P-g_m)$ is a  projective pure
submodule of the projective left $R$-module $\mathrm{Hom}_R(P,R)$.

Notice that, for any $m$, $\overline{g_m}=\alpha \circ \beta$ and
$\overline{\mathrm{Id}_P-g_m}=\mathrm{Id}_{P/PI}-\alpha\circ \beta$.
Therefore, $P'_n/P'_nI\cong Q/QI$, hence $P'_n\cong Q$, and
$$Q'_n/IQ'_n\cong \mathrm{Hom}_{R/RI}((\mathrm{Id}_{P/PI}-\alpha\circ
\beta)P/PI,R/I).$$ In particular, the isomorphism classes of $P'_n$
and $Q'_n$, respectively, do not depend on $n$.
\end{Rem}

Combining Proposition \ref{sequence} with Remark \ref{left} we
obtain the following theorem which is a slight refinement of
\cite[Theorem 7.1]{plans}.

\begin{Th}\label{flatproj} Let $R$ be a ring, let $P$ be a finitely generated projective right $R$-module, and let $I\subseteq J(R)$ be a two-sided ideal of $R$.
Assume that there is a split exact sequence of right $R/I$ modules
\[0\to X\to P/PI\to X'\to 0.\]
Then the following statements are equivalent,
\begin{itemize}
\item[(i)] There exists a (countably generated) projective right
$R$-module $Q$ such that $Q/QI\cong X$.
\item[(ii)] There exists a (countably generated) projective left
$R$-module $Q'$ such that $Q'/IQ'\cong \mathrm{Hom}_{R/I}(X',R/I)$.
\end{itemize}
When the above equivalent statement  hold $Q$ is isomorphic to a
pure submodule of $P$, and $Q'$ is isomorphic to a pure submodule
of $\mathrm{Hom}_R(P,R)$. Moreover, $Q$ is finitely generated if and
only if $Q'$ is finitely generated if and only if there exists a
projective right $R$-module $P'$ such that $P'/P'I\cong X'$.
\end{Th}

Now we are going to state some of the results above in terms of
monoids of projectives. More precisely, in terms of pre-ordered
monoids of projectives.

We recall that over a commutative monoid $M$ there is a pre-order
relation called \emph{the algebraic preorder} on $M$ defined by
$x\ge y$, for $x$, $y\in M$, if and only if $x=y+z$ for some $z\in
M$.

For example, over $(\No^*)^k$ the algebraic order is the
component-wise order, which is even a partial order. When the monoid
is $V^*(R)$ for some ring $R$, $\langle Q\rangle \le \langle
P\rangle$ if and only if $Q$ is isomorphic to a direct summand of
$P$.

In terms of monoids of projective modules Corollary \ref{full}
essentially says that for elements in $V(R)$ the algebraic preorder
is respected modulo $J(R)$. We state this in a precise way in the
next result.

\begin{Cor} \label{fullmonoid} Let $R$ be a ring, and let $I$ be a two-sided ideal of
$R$ contained in $J(R)$. Let $\pi \colon R\to R/I$ denote the
projection, and let $\stackrel{\sim}{\pi}\colon V^*(R)\to V^*(R/I)$
denote the induced homomorphism of monoids. If $x\in V^*(R)$, $y\in
V(R)$ are such that there exist $c\in V^*(R/I)$ satisfying that
$\stackrel{\sim}{\pi}(x)=\stackrel{\sim}{\pi}(y)+c$ then there
exists $z\in V^*(R)$ such that $\stackrel{\sim}{\pi}(z)=c$ and
$x=y+z$.
\end{Cor}

In  general, for a semilocal ring $R$,   the monoid $V^*(R)$ is
isomorphic to a submonoid of $(\No^*)^k$. In view of Theorem
\ref{Pavel}, the algebraic order of $(\No^*)^k$ induces an order on
$V^*(R)$ that is translated in terms of projective modules over $R$
by $\langle Q\rangle \le \langle P\rangle$ if and only if there
exists a pure monomorphism $f\colon Q\to P$ if and only if $Q/QJ(R)$
is a direct summand of $P/PJ(R)$. By \cite{P1}, the relation $\le$
is antisymmetric. This partial order relation defined on $V^*(R)$
restricts to the usual algebraic order over $V(R)$, but not on
$W(R)$ when $V(R)\subsetneq W(R)$.

\begin{Cor} \label{wrvr} Let $R$ be a semilocal ring, fix $\varphi\colon R\to S$ an onto ring homomorphism to a semisimple artinian ring $S$
such that $\mathrm{Ker}\, \varphi =J(R)$. Then
\begin{itemize}
\item[(i)] $x\in W(R)\setminus V(R)$ if and only if $x$ is incomparable
(with respect to the algebraic order) with $n\langle R\rangle$ for
any $n\ge 1$ if and only if there exist $n\ge 1$ such that $n\cdot
\dim _\varphi\langle R\rangle-\dim _\varphi (x)\in \dim_\varphi
W({}_RR)\setminus \dim_\varphi V(R)$.

\item[(ii)] $V(R)=W(R)\cap W({}_RR).$
\end{itemize}
\end{Cor}

\begin{Proof} Since over a semisimple artinian ring any exact sequence splits, the statement follows by applying Theorem
\ref{flatproj}.
\end{Proof}

\begin{Rem} \label{finiteandinequalities} Corollary \ref{wrvr} implies that,   if $\dim _\varphi
V^*(R)\subseteq (\No ^*)^k$ is a monoid defined by inequalities and
$$\dim_\varphi (\langle R\rangle)\in \dim _\varphi V^*(R)\cap D(\dim
_\varphi V^*(R)),$$ the elements of the semigroup $\dim _\varphi
W(R)\setminus \dim _\varphi V(R)$ must be the elements of $\No^k$
such that some of the inequalities they satisfy are strict. So that
$$\dim _\varphi V(R)=\dim _\varphi V^*(R)\cap D(\dim _\varphi
V^*(R))\cap \No^k=\dim _\varphi W(R)\cap \dim _\varphi
W({}_RR).$$
\end{Rem}

In terms of order relations on the monoids we have the following
Corollary.

\begin{Cor}\label{algorder} Let $R$ be a   semilocal ring.
Consider the following relation over $V^*(R)$, $\langle P\rangle \le
\langle Q\rangle $ if and only if $P/PJ(R)$ is isomorphic to a
direct summand of $Q/QJ(R)$. Then
\begin{itemize}
\item[(i)] $\langle P\rangle \le
\langle Q\rangle $ if and only if there exists a pure embedding
$f\colon P\to Q$.
\item[(ii)] $\le $ is a partial order relation that refines the algebraic order on $V^*(R)$.
\item[(iii)] If, in addition, $R$ is noetherian then the partial
order induced  by $\le $ over $V^*(R)$ is the algebraic order.
\end{itemize}
\end{Cor}

\begin{Proof} $(i).$ If  $\langle P\rangle \le
\langle Q\rangle $ then there exists a splitting monomorphism
$\overline{f}\colon P/PJ(R)\to Q/QJ(R)$ which by
Theorem~\ref{Pavel}(i) lifts to a pure monomorphism $f\colon P\to
Q$.  Conversely, if $f\colon P\to Q$ is a pure monomorphism of right
$R$-modules then the induced map $f\otimes _R R/J(R)\colon P\otimes
_R R/J(R)\to Q \otimes _R R/J(R)$ is a pure monomorphism of
$R/J(R)$-modules. Since $R/J(R)$ is semisimple, $f\otimes _RR/J(R)$
is a split monomorphism.

$(ii).$ It is clear that $\le $ is reflexive and transitive. As it
is already observed in $\cite{P1}$, Theorem~\ref{Pavel} implies that
$\le$  is also antisymmetric.

If $P$ is isomorphic to a direct summand of $Q$, then $P/PJ(R)$ is
also isomorphic to a direct summand of $Q/QJ(R)$. Hence $\langle
P\rangle \le \langle Q\rangle$, that is, $\le$ refines the algebraic
order on $V^*(R)$.

$(iii).$ It is a consequence of the realization Theorem~\ref{main}.
\end{Proof}

We shall see in Examples \ref{noniso} that the monoid $V^*(R)$ does
not determine $V^*({}_RR)$. Theorem \ref{flatproj}, or \cite[Theorem
7.1]{plans}, combined with Theorem \ref{Pavel}(ii) implies that for
a semilocal ring $W(R)$ does determine  $W({}_RR)$.

\begin{Cor}\label{uniquewr} For $i=1,2$, let $R_i$ be a semilocal
ring and let $\varphi _i\colon R_i\to M_{n_1}(D^i_1)\times \cdots
\times M_{n_k}(D^i_k)$ be an onto ring homomorphism such that
$\mathrm{Ker}\, \varphi _i=J(R_i)$ and $D^i_1,\dots, D^i_k$ are
division rings.

Then $\dim _{\varphi _1} W(R_1)=\dim _{\varphi _2} W(R_2)$ if and
only if $\dim _{\varphi _1} W({}_{R_1}R_1)=\dim _{\varphi _2}
W({}_{R_2}R_2)$.
\end{Cor}

\begin{Proof} By symmetry, it is enough to prove that if $\dim _{\varphi _1} W(R_1)=\dim _{\varphi _2} W(R_2)$ then  $\dim
_{\varphi _1} W({}_{R_1}R_1)\subseteq \dim _{\varphi _2}
W({}_{R_2}R_2)$.

Let $x\in \dim _{\varphi _1} W({}_{R_1}R_1)$. There exists $m\in \N$
such that $x\le m(n_1,\dots ,n_k)$. By Theorem \ref{flatproj},
$y=m(n_1,\dots ,n_k)-x \in \dim _{\varphi _1} W(R_1)=\dim _{\varphi
_2} W(R_2)$. Applying again Theorem \ref{flatproj}, we deduce that
$x=m(n_1,\dots ,n_k)-y\in \dim _{\varphi _2} W({}_{R_2}R_2)$.
\end{Proof}

\section{Some   examples} \label{examples}

Gerasimov and Sakhaev gave the first example of a semilocal ring
such that $V(R)\varsubsetneqq W(R)$. The final step for the
computation of $V^*(R)$ was made in \cite{DPP}. We want to start
this section stating the main properties of this example as it  is
one of the basic tools to prove our realization Theorem
\ref{pbinequalities}.

\begin{Th}\label{GS} \emph{(\cite{GS}, \cite{DPP})} Let $F$ be any field. There exists a semilocal $F$-algebra $R$ with an onto ring morphism  $\varphi\colon R\to F\times
F$ with $\mathrm{Ker}\, \varphi =J(R)$ and such that all finitely
generated projective modules are free but
\[\begin{array}{rcl}\dim _\varphi
W(R_R)&=&\{ (x,y)\in \No\mid x\ge y\}=(1,1)\No+(1,0)\No\\\dim
_\varphi V^*(R_R)&=&\left(\dim _\varphi W(R_R)\right)\No^*=\{
(x,y)\in \No^*\mid x\ge y\}
\end{array}\] and
\[\begin{array}{rcl}\dim _\varphi W({}_RR)&=&\{ (x,y)\in \No\mid y\ge x\}=(1,1)\No+(0,1)\No \\ \dim
_\varphi V^*({}_RR)&=&\left(\dim _\varphi W({}_RR)\right)\No^*=\{
(x,y)\in \No^*\mid y\ge x\}\end{array}\] In particular, any
projective module over $R$ is a direct sum of indecomposable
projective modules that are finitely generated modulo $J(R)$.
\end{Th}

It is quite  an interesting question to determine the structure of
$V^*(R)$ for a general semilocal ring. But right now it seems to be
too challenging  even for   semilocal rings  $R$ such that
$R/J(R)\cong D_1\times D_2$ where  $D_1$, $D_2$ are division rings.
Now we provide some examples of such rings to illustrate
Theorem~\ref{pbinequalities} and the difficulties that  appear in
the general case. We first observe that, since $k=2$ and $(1,1)$ is
the order unit of $\dim _\varphi V(R)$, to have some room for
interesting behavior of countably generated projective modules all
finitely generated projective modules must be free.

\begin{Lemma}\label{free} Let $R$ be a semilocal ring such that $R/J(R)\cong D_1\times D_2$ for
suitable division rings $D_1$ and $D_2$. Fix $\varphi \colon R\to
D_1\times D_2$ an onto ring homomorphism such that $\mathrm{Ker}\,
\varphi =J(R)$. If $R$ has non-free finitely generated projective
right (or left) modules then there exists $n\in \N$ such that $\dim
_\varphi V(R)$ is the submonoid of $\No ^2$ generated by $(1,1)$,
$(n,0)$ and $(0,n)$. In this case,
\[\dim _\varphi V^*(R)=(1,1)\No^*+(n,0)\No^*+(0,n)\No^*=\{(x,y)\in \No ^*\mid x+(n-1)y\in n\No ^*\}.\]
Therefore, all projective modules are direct sum of finitely
generated projective modules.
\end{Lemma}

\begin{Proof} Note that $\dim _\varphi (\langle R\rangle )=(1,1)$. So that $(1,1)\in A = \dim _\varphi
V(R)$.

Let $P$ be a non-free finitely generated projective right
$R$-module, and let $\dim _\varphi (\langle P \rangle)=(x,y)$. As
$P$ is not free, either $x>y$ or $x<y$. Assume $x>y$, then
\[(x,y)=(x-y,0)+y(1,1)\in A\qquad (*).\] Since, by
Corollary~\ref{full} or its monoid version Corollary
\ref{fullmonoid}, $A$ is a full affine submonoid of $\No ^2$ we
deduce that $(x-y,0)\in A$ and also that
$(0,x-y)=(x-y)(1,1)-(x-y,0)\in A$. If $x<y$ we deduce, in a
symmetric way that $(y-x,0)$ and $(0,y-x)$ are elements of $A$.

Choose $n\in \N$ minimal with respect to the property $(n,0)\in A$,
and note that then also $(0,n)\in A$. We claim that
\[A=(1,1)\No+(n,0)\No +(0,n)\No.\]
We only need to prove that if $(x,y)\in A$ then it can be written as
a linear combination, with coefficients in $\No$ of $(1,1)$, $(n,0)$
and $(0,n)$. In view of the previous argument, it suffices to show
that if $(x,0)\in A$ then $(x,0)\in (n,0)\No$. By the division
algorithm $(x,0)=(n,0)q+(r,0)$ with $q\in \No$ and $0\le r<n$. As
$A$ is a full affine submonoid of $\No^2$ we deduce that $(r,0)\in
A$. By the choosing of $n$, $r=0$ as desired.

Let $P_1$ be a finitely generated right $R$-module such that $\dim
_\varphi (\langle P_1\rangle)=(n,0)$, and let $P_2$ be a finitely
generated right $R$-module such that $\dim _\varphi (\langle
P_2\rangle)=(0,n)$.

Let $Q$ be a countably generated projective right $R$-module that is
not finitely generated.  Let $\dim _\varphi (\langle Q\rangle
)=(x,y)\in \No ^*$. We want to show that
\[(x,y)\in (1,1)\No^*+(n,0)\No^*+(0,n)\No^*\]

If $x=y$ then $(x,y)=x(1,1)$ and, by Theorem~\ref{Pavel}(ii), $Q$ is
free. If $x>y$ then $y\in \No$ and $(x,y)=(x-y,0)+y(1,1)$, combining
Theorem~\ref{Pavel}(ii) with Lemma~\ref{full} we deduce that
$Q=yR\oplus Q'$ with $Q'$ such that $\dim _\varphi (\langle
Q'\rangle )=(z,0)$ where $z=x-y$. If $z<\infty$ then, by
Theorem~\ref{Pavel}(ii), $nQ'\cong zP_1$ hence $Q'$, and $Q$, are
finitely generated. If $z=\infty$, by Theorem~\ref{Pavel}(ii),
$Q'\cong P_1^{(\omega)}$. Hence $(x,y)=\infty \cdot (n,0)+y(1,1)$.
The case $x<y$ is done in a symmetric way.

It is not difficult to check that the elements of $\dim _\varphi
V^*(R)$ are the solutions in $\No^*$ of $x+(n-1)y\in n\No ^*$.
\end{Proof}

Now we will list all the possibilities for the monoid $V^*(R)$
viewed as a submonoid of $V^*(R/J(R))$ when $R$ is a noetherian ring
such that $R/J(R)\cong D_1\times D_2$, for $D_1$ and $D_2$ division
rings, and all finitely generated projective modules are free. In
view of Theorem \ref{main} this is equivalent to classify the
submonoids of $(\No^*)^2$ containing $(1,1)$ and that are defined by
a system of equations. Though the presentation of the monoid as
solutions of equations is quite attractive there is an alternative
one that, even being technical, is more useful to work with.

\begin{Def} \label{defsystems}
Fix $k \in \mathbb{N}$ and an order unit $(n_1,\dots,n_k) \in
\mathbb{N}^k$. A \emph{system of supports}
$\mathcal{S}(n_1,\dots,n_k)$ consists of a collection $\cal S$ of
subsets of $\{1,\dots,k\}$ together with a family of commutative
monoids $\{A_I, I \in \cal S\}$ such that the following conditions
hold
\begin{enumerate}
\item[(i)] $\emptyset$ and $\{ 1,\dots,k\}$ are elements of $\cal S$.
\item[(ii)] For any $I \in \cal S$, $A_I$ is a  submonoid of
$\mathbb{N}_0^{\{1,\dots,k\} \setminus I}$. The monoid
$A_{\{1,\dots,k\}}$ is the trivial monoid and $(n_1,\dots,n_k) \in
A_{\emptyset}$.
\item[(iii)]  $\mathcal{S}$ is closed under unions, and if $x\in A_I$ for some $I\in \mathcal{S}$ then
$I\cup \supp(x)\in \cal S$. In particular $\{1,\dots ,k\}\in
\mathcal{S}$.
\item[(iv)] Suppose that $I,K \in \cal S$ are  such that $I \subseteq K$ and
let $p \colon \mathbb{N}_0^{\{1,\dots,k\}\setminus I} \to
\mathbb{N}_0^{\{1,\dots,k\} \setminus K}$ be the canonical
projection. Then $p(A_I) \subseteq A_K$.
\end{enumerate}

If in addition, for any $I\in \Scal$, the submonoids $A_I$ are full
affine submonoids of $\No ^{\{1,\dots,k\} \setminus I}$ then
$\mathcal{S}(n_1,\dots,n_k)$ is said to be a \emph{full affine
system of supports}.
\end{Def}

\begin{Rem} \label{remsupports} Given a system of supports $\mathcal{S}(n_1,\dots,n_k)=\{A_I, I \in \cal S\}$ we can
associate to it a monoid.  Consider the subset $M(\mathcal{S})$ of
$(\No ^*)^k$ defined by $x\in M(\mathcal{S})$ if and only if $I =
\infsupp (x) \in \cal S$ and $p_I(x) \in A_I$, where if
$x=(x_1,\dots ,x_k)$ then
$$\infsupp (x)=\{i\in \{1,\dots,k\}\mid x_i=\infty\},$$
and $p_I\colon\No^k\to \mathbb{N}_0^{\{1,\dots,k\} \setminus I}$
denotes the canonical projection.

By \cite[Theorem 7.7]{HP}, $\mathcal{S}(n_1,\dots,n_k)$ is a full
affine system of supports if and only if $M(\mathcal{S})$ is a
monoid defined by equations and containing $(n_1,\dots ,n_k)$.
\end{Rem}

We recall that a module is superdecomposable if it has no
indecomposable direct summand. By Theorem~\ref{main} and
Lemma~\ref{full}, in our context superdecomposable modules are
relatively frequent as they correspond to the elements $x\in
M\subseteq (\No^*)^k$ such that, for any $y\in M\cap \No^k$,
$\supp\, (y)\nsubseteqq \supp\, (x)$.

\begin{Ex} \label{nk2} Let $R$ be a semilocal noetherian ring such that there exists   $\varphi\colon R\to D_1\times D_2$, an
onto ring morphism with $\mathrm{Ker}\, \varphi=J(R)$, where  $D_1$
and $D_2$ are division rings. Assume that  all finitely generated
projective right $R$-modules are free. Hence $\dim _\varphi V(R)=
(1,1)\No$, and its order unit is $(1,1)$. Then there are the
following possibilities for $\dim_\varphi V^*(R)$:
\begin{itemize}
\item[(0)] All projective modules are free, so that
$M_0=\dim _\varphi V^*(R)=(1,1)\No^*$. Note that $M_0$ is the set of
solutions   $(x,y)\in (\No^*) ^2$ of the equation $x=y$.
\item[(1)] $M_1=\dim _\varphi V^*(R)=(1,1)\No^*+ (0,\infty)\No^*$. So that, $M_1$ is the set of solutions $(x,y)\in (\No^*)^2$
of the equation $x+y=2y$.

Note that for such an $R$ there exists a countably generated
superdecomposable projective right $R$-module $P$ such that $\dim
_\varphi (\langle P\rangle )=(0,\infty)$. Then any countably
generated projective right $R$ module $Q$ is isomorphic to
$R^{(n)}\oplus P^{(m)}$ for suitable $n\in \No ^*$ and $m\in
\{0,1\}$.

\item[(1')] $M_1'=\dim _\varphi V^*(R)=(1,1)\No^*+ (\infty,0)\No^*$. So that, $M_1'$ is the set of solutions $(x,y)\in (\No^*)^2$
of the equation $x+y=2x$.

\item[(2)] $M_2=\dim _\varphi V^*(R)=(1,1)\No^*+ (\infty,0)\No^*+(0,\infty)\No^*$.
So that, $M_2$ is the set of solutions $(x,y)\in (\No^*)^2$ of the
equation $2x+y=x+2y$.

Note that for such an $R$ there exist two  countably generated
superdecomposable projective right $R$-modules $P_1$ and $P_2$ such
that $\dim _\varphi (\langle P_1\rangle )=(0,\infty)$ and $\dim
_\varphi (\langle P_2\rangle )=(\infty ,0)$. Any countably generated
projective right $R$ module $Q$ satisfies that there exist $n\in
\No$ and $m_1,m_2\in \{0,1\}$ such that $Q=R^{(n)}\oplus
P_1^{(m_1)}\oplus P_2^{(m_2)}$.
\end{itemize}
\end{Ex}

\begin{Proof} In view of Theorem~\ref{main} and Remark~\ref{remsupports} we must describe all the possibilities for full
affine systems of supports of $\{1,2\}$ such that $A_{\emptyset} =
(1,1) \mathbb{N}_0$. Since the set of supports of a system of
supports at least contains $\emptyset$ and $\{1,2\}$ there are just
four possibilities.

Since the image of the projections of $A_\emptyset$  on the first
and on  the second component is $\N_0$, all the monoids $A_I$ in the
definition of system of supports are determined by $A_\emptyset$.

Case $(0)$ is the one in which $M_0=A_\emptyset+\infty\cdot
A_\emptyset$. According to Remark \ref{allfg} (3), in this case all
projective modules are direct sum of finitely generated
(indecomposable) modules.

In cases $(1)$ and $(1')$ there are $3$ different  supports for the
elements in the monoid, and in case $(2)$ there are $4$.
\end{Proof}

Now we give some examples whose existence is a direct consequence of
Theorem \ref{pbinequalities}.

\begin{Ex}\label{noniso} Let $F$ be any field. In all the statements $R$ denotes a semilocal F-algebra, and
  $\varphi\colon R\to E\times E$ denotes an onto ring homomorphism
such that $\mathrm{Ker}\, \varphi=J(R)$ and $E$ is a suitable field
extension of $F$. Fix $n\in \N$. Then there exist $R$ and $\varphi$
such that
\begin{itemize}
\item[(i)]
\[N=\dim_\varphi V^*(R_R)=(1,1)\No^*+(n,0)\No^*=\{(x,y)\in (\No^*)^2\mid x\ge y \mbox{ and }x+(n-1)y\in n\No ^*\}\]
\[D(N)=\dim_\varphi V^*({}_RR)=(1,1)\No^*+(0,n)\No^*=\{(x,y)\in (\No^*)^2\mid x\le y \mbox{ and }x+(n-1)y\in n\No ^*\}\]
For $n=1$, we recover the situation in \cite{GS}. Note that over $R$
all projective modules are direct sum of indecomposable projective
modules.
\item[(ii)]
\[\dim_\varphi V^*(R_R)=N+(0,\infty)\No^*=\{(x,y)\in (\No^*)^2\mid 2x+y\ge 2y+x \mbox{ and }x+(n-1)y\in n\No ^*\}\]
\[\dim_\varphi V^*({}_RR)=D(N)+(\infty ,0)\No^*=\{(x,y)\in (\No^*)^2\mid 2x+y\le 2y+x \mbox{ and }x+(n-1)y\in n\No ^*\}\]
In this case $R$ has a superdecomposable projective right $R$-module
and a superdecomposable projective left $R$-module.
\item[(iii)]
\[\dim_\varphi V^*(R_R)=N+(0,\infty)\No^*=\{(x,y)\in (\No^*)^2\mid x+y\ge 2y \mbox{ and }x+(n-1)y\in n\No ^*\}\]
\[\dim_\varphi V^*({}_RR)=D(N)=\{(x,y)\in (\No^*)^2\mid x+y\le 2y \mbox{ and }x+(n-1)y\in n\No ^*\}\]
In this situation $R$ has a superdecomposable projective right
$R$-modules but every projective left $R$-module is a direct sum of
indecomposable modules.
\item[(iv)]
\[\dim_\varphi V^*(R_R)=(1,1)\No^*+(\infty,0)\No^*=\{(x,y)\in (\No^*)^2\mid 2x =  x+y\mbox{ and }x\ge y\}\] and
\[\dim_\varphi V^*({}_RR)=(1,1)\No^*=\{(x,y)\in (\No^*)^2\mid 2x= x+y \mbox{ and }x\le y\}.\] Therefore, all projective left
$R$-modules are free hence they are a direct sum of finitely
generated modules but this is not true for projective right
$R$-modules. In particular, $V^*(R_R)$ and $V^*({}_RR)$ are not
isomorphic.
\end{itemize}
In the first three examples $V(R)\varsubsetneq
W(R)=(1,1)\No+(n,0)\No\cong W({}_RR)$. In the fourth example, as
Theorem \ref{flatproj} implies, $V(R)=W(R)=W({}_RR)$.
\end{Ex}

\begin{Proof} After Theorem \ref{pbinequalities} what is left to do is to check the generating sets of the monoids. But all the computations are straightforward.

In $(iv)$ to prove that $V^*(R)$ is not isomorphic to $V^*({}_RR)$
just count the number of idempotent elements in both monoids.
\end{Proof}

\begin{Rem}
Examples \ref{noniso}(ii) and (iii) answer a problem mentioned in
\cite[page~3261]{DPP}, and Example \ref{noniso}(iv) answers a
problem in \cite[page~310]{FS}.

Following the notation of Examples \ref{noniso} and under the same
hypothesis, the first place where it was shown that there could be a
non finitely generated projective module $P$ such that
$\mathrm{dim}_\varphi (\langle P\rangle)= (n,0)$ for a given $n>1$
was in \cite{Sa}.

The monoid $M=N+(0,\infty)\No^*$ is described in Examples
\ref{noniso}(ii) and (iii) in two different ways as a monoid given
by a system of inequalities. Both descriptions result in different
monoids $D(M)$.
\end{Rem}

Now we give an example such that $W(R)\not \cong W({}_RR)$ and
$V^*(R)\not \cong V^*({}_RR)$. It also shows that Corollary
\ref{full} fails also for the semigroup $W(R)\setminus V(R)$, so
that in Theorem \ref{flatproj} we cannot just assume that $P$ is
finitely generated modulo the Jacobson radical.

\begin{Ex}\label{nondivisibility} Fix $1\le n\in \N$. Let $F$ be any field. There exist a   semilocal F-algebra $R$, a suitable field extension $E$ of $F$ and
    an onto ring homomorphism $\varphi\colon R\to E\times M_n(E)$
such that $\mathrm{Ker}\, \varphi=J(R)$ and
\[\dim_\varphi
V^*(R)=(1,n)\No^*+\cdots+(1,0)\No^*=\{(x,y)\in (\No^*)^2\mid nx\ge y
\}\]
\[\dim_\varphi V^*({}_RR)=(1,n)\No^*+(0,1)\No^*=\{(x,y)\in (\No^*)^2\mid nx\le y \}.\]
Therefore $W(R)=(1,n)\No+\cdots+(1,0)\No$ and
$W({}_RR)=(1,n)\No+(0,1)\No$ which are non isomorphic monoids
provided $n\ge 2$.

Notice that the $(1,0),\dots ,(1,n-1)$ are minimal elements of
$W(R)$ and of $W(R)\setminus V(R)$ so that they are incomparable.
\end{Ex}

\begin{Proof} The existence of the semilocal ring follows from
Theorem \ref{pbinequalities}. We show that the two  monoids have the
required set of generators.

Let $M=\{(x,y)\in (\No^*)^2\mid nx\ge y \}$. It is clear that
$(1,n)\No^*+\cdots+(1,0)\No^*\subseteq M$. If $(x,y)\in M\cap \No^k$
then $y=n\cdot k+y'$ for some $k,y'\in \No$ and $0\le y'<n$.
Therefore, if $x=k$, $(x,y)=k(1,n)$. If $x>k$ then
$(x,y)=k(1,n)+(x-k-1)(1,0)+(1,y')$ provided $y'>0$, otherwise
$(x,y)=k(1,n)+(x-k)(1,0)$. In the three cases we conclude that
$(x,y)\in (1,n)\No+\cdots+(1,0)\No$. For elements with nonempty
infinite support the inclusion is clear.

If $(x,y) \in D(M)\cap \No^k$ then $(x,y)=x(1,n)+(y-nx)(0,1)$ which
proves that $D(M)=(1,n)\No^*+(0,1)\No^*$.

The monoids $W(R)$ and $W({}_RR)$ have the same number of minimal
elements if and only if $n=1$. Therefore they cannot be isomorphic
for $n\ge 2$.
\end{Proof}

\section{Monoids defined by inequalities}\label{monoids}

We think on $(\No^*)^k$ and of $\No^k$ as ordered monoids with the
order relation given by the algebraic order. That is,
$(x_1,\dots,x_k)\le (y_1,\dots,y_k)$ if and only if $x_i\le y_i$ for
any $i=1,\dots ,k$.

We recall that a monoid $M$ is said to be unperforated if, for every
$n\in \N$, it satisfies the following properties:
\begin{itemize}
\item[(1)] For any $x$, $y\in M$, $nx\le ny$ implies  $x\le y$;
\item[(2)] for any $x$, $y\in M$, $nx=ny$ implies  $x= y$.
\end{itemize}
where $\le$ denotes the algebraic preordering on $M$.

\begin{Prop}\label{characterizations} (\cite[Proposition~2]{KL}) Let  $A$  be a commutative cancellative monoid such that $U(A)=\{0\}$. Then
the following statements are equivalent;
\begin{itemize}
\item[(i)] $A$ is finitely generated and unperforated.
\item[(ii)] There exist $k\ge 1$, a monoid embedding $f\colon A\to \No^k$ and $E_1,$ $E_2\in M_{\ell \times
k}(\No)$ such that $f(A)$ is the set of solutions in $\No^k$ of the
system
\[E_1\left(\begin{array}{c}t_1\\\vdots \\
t_k\end{array}\right)= E_2\left(\begin{array}{c}t_1\\\vdots \\
t_k\end{array}\right)\]
\item[(iii)] There exist $m\ge 1$, a monoid embedding $g\colon A\to \No^m$, $D\in M_{n\times m}(\No)$,
$E_1,\, E_2\in M_{\ell \times m}(\No)$  and $m_1,\dots ,m_n\in \N$ ,
$m_i\ge 2$ for any $i\in \{1,\dots ,n\}$, such that $g(A)$ is the
set of solutions in $\No^m$ of the system
\[D\left(\begin{array}{c}t_1\\\vdots \\ t_m\end{array}\right)\in \left(\begin{array}{c}m_1\No\\
\vdots \\ m_n\No\end{array} \right)  \qquad \mbox{and} \qquad
E_1\left(\begin{array}{c}t_1\\\vdots \\
t_m\end{array}\right)= E_2\left(\begin{array}{c}t_1\\\vdots \\
t_m\end{array}\right) \]
\item[(iv)] There exist $s\ge 1$, a monoid embedding $h\colon A\to \No^s$, $D\in M_{n\times s}(\No)$,  $E_1,\, E_2\in M_{\ell \times
s}(\No)$  and $m_1,\dots ,m_n\in \N$ , $m_i\ge 2$ for any $i\in
\{1,\dots ,n\}$, such that $h(A)$ is the set of solutions in $\No^s$
of the system
\[D\left(\begin{array}{c}t_1\\\vdots \\ t_s\end{array}\right)\in \left(\begin{array}{c}m_1\No\\
\vdots \\ m_n\No\end{array} \right)  \qquad \mbox{and} \qquad
E_1\left(\begin{array}{c}t_1\\\vdots \\
t_s\end{array}\right)\ge E_2\left(\begin{array}{c}t_1\\\vdots \\
t_s\end{array}\right) \]
\end{itemize}
\end{Prop}

\begin{Proof} For further quoting we give the proof of the equivalence of $(iii)$ and $(iv)$. It is clear that the monoids in $(iii)$ can be described
as the set of solutions of a system of congruences and inequalities
as the ones appearing in $(iv)$.

Conversely, let $A$ be a submonoid of $\No ^s$ that is the set of
solutions in $\No ^s$ of  the system in $(iv)$. Consider the monoid
morphism $g\colon A\to \No ^{s+\ell}$ defined by
$$g(a_1,\dots ,a_s)=\left( a_1,\dots a_s,\sum _{i=1}^se^1_{1i}a_i-\sum
_{i=1}^se^2_{1i}a_i,\dots ,\sum _{i=1}^se^1_{\ell i}a_i-\sum
_{i=1}^se^2_{\ell i}a_i\right)$$ where $(a_1,\dots ,a_s)\in A$ and,
for $k=1,2$, $e_{ij}^k$ denotes the $i$-$j$-entry of the matrix
$E_k$.

Then $g(A)$ is the set of solutions in $\No^{s+l}$ of the system
\[D\left(\begin{array}{c}t_1\\\vdots \\ t_s\end{array}\right)\in \left(\begin{array}{c}m_1\No\\
\vdots \\ m_n\No\end{array} \right)  \qquad \mbox{and} \qquad
E_1\left(\begin{array}{c}t_1\\\vdots \\
t_s\end{array}\right)=E_2 \left(\begin{array}{c}t_1\\\vdots \\
t_s\end{array}\right)+ \left(\begin{array}{c}t_{s+1}\\\vdots \\
t_{s+\ell}\end{array}\right). \] So that $A$ is also a monoid of the
type appearing in $(iii)$.
\end{Proof}

The embeddings of $(iii)$ are the full affine embeddings. We recall
that if $g(A)$ has an order unit $(n_1,\dots, n_m)$ of $\N ^m$ then
$g(A)$ can be realized as $\dim _\varphi (V(R))$ for some semilocal
ring $R$ such that $R/J(R)\cong M_{n_1}(D_1)\times \cdots \times
M_{n_m}(D_m)$ for suitable division rings $D_1,\dots ,D_m$
\cite{FH}.

We stress that not all finitely generated submonoids of $\No^k$ are
unperforated. Consider, for example, $N=(1,1)\No+(2,0)\No+(3,0)\No$.
In $N$, $2(2,0)\le 2(3,0)$ but $(2,0)$ and $(3,0)$ are incomparable
in $N$.

In the next lemma we study  monoids defined by a system of equations
and monoids defined by a system of inequalities.

\begin{Lemma} \label{monoequality} Let   $M$ be a submonoid of $(\No^*)^k$ defined by a system of
inequalities
\[D\cdot T\in \left(\begin{array}{c}m_1\No^*\\
\vdots \\ m_n\No^*\end{array} \right) \quad(*) \qquad \mbox{and}
\qquad E_1\cdot T\ge E_2\cdot T \quad(**)\] where $T=(t_1,\dots
,t_k)^t$, $D\in M_{n\times k}(\No)$,  $E_1,\, E_2\in M_{\ell \times
k}(\No)$  and $m_1,\dots ,m_n\in \N$ , $m_i\ge 2$ for any $i\in
\{1,\dots ,n\}$. Let $A$ be the submonoid of $M$ whose elements are
the solutions in $\No^k$ of
\[D\cdot T\in \left(\begin{array}{c}m_1\No\\
\vdots \\ m_n\No\end{array} \right) \quad  \qquad \mbox{and} \qquad
E_1\cdot T= E_2\cdot T\] Then,
\begin{itemize}
\item[(i)] $M$ and $D(M)$ are finitely generated monoids.
\item[(ii)]$A=M\cap D(M)\cap\No^k$.
\item[(iii)] For any $m\in M$ and $a\in A$, if  there exists $m'\in(\No^*)^k$ such that  $m=a+m' $
then $m' \in M$.
\end{itemize}
\end{Lemma}

\begin{Proof} $(i)$ Consider the monoid $N$ defined the system of
equations
\[D'\cdot T'\in \left(\begin{array}{c}m_1\No^*\\
\vdots \\ m_n\No^*\end{array} \right) \quad(*) \qquad \mbox{and}
\qquad E_1\cdot T= E_2\cdot T
+\left(\begin{array}{c}t_{k+1}\\\vdots\\t_{k+\ell}
\end{array}\right)\quad(**)\] where $T'=(t_1,\dots ,t_k, t_{k+1},\cdots ,t_{k+\ell})^t$ and $D'=(D\vert0)\in M_{n\times (k+\ell)}(\No)$. By
\cite[Example 7.6]{HP}, $N$ is a finitely generated monoid.

Let $p\colon (\No^*)^{k+\ell}\to (\No^*)^k$ denote the projection
onto the first $k$ components. It is easy to see that $p(N)=M$, so
that  $M$ is finitely generated.

Statements $(ii)$ and $(iii)$ are clear.
\end{Proof}

In contrast with Proposition \ref{characterizations}, the monoid $N$
appearing in the proof of Lemma \ref{monoequality} need not be
isomorphic to $M$.

In general, as the following basic example shows, a monoid defined
by inequalities may not be isomorphic to a monoid defined by a
system of equations. Therefore the equivalence of statements (ii),
(iii) and (iv) in Proposition \ref{characterizations} does not
extend to submonoids on $(\No^*)^k$.

\begin{Ex} Let $M$ be the submonoid of $(\No^*)^2$ that is the set
of  solutions of $x\ge y$. Then $M$ is not isomorphic to a monoid
defined by a system of equations.
\end{Ex}

\begin{Proof} In order to be able to manipulate this monoid we need
to think on the language of system of supports, see Definition
\ref{defsystems} and Remark~\ref{remsupports}.

First note that $M=(1,1)\No+(1,0)\No+(\infty
,0)\No+(\infty,\infty)\No$. The elements $c=(\infty ,0)$ and
$d=(\infty,\infty)$ are nonzero elements satisfying that $2c=c$,
$2d=d$ and $d+c=d$. Therefore, if $h\colon M\to N$ is a monoid
morphism and $N$ is a monoid defined by a system of equations,
$h(c)$ and $h(d)$ must be elements such that its support coincides
with its infinite support and, moreover, $\supp h(c)\subseteq \supp
h(d)$. If $h$ is bijective, then $h(c)$ and $h(d)$ are the only
non-zero elements of $N$ such that its support coincides with its
infinite support. So that if we think on the presentation of $N$ as
a system $\mathcal{S}$ of supports, we deduce that there are only
three different sets in $\mathcal{S}$, that is $\emptyset$ , $\supp
h(c)$ and $\supp h(d)$. Moreover, $\supp h(c)\subsetneq \supp h(d)$

On the other hand, since $(1,0)+c=c$, we deduce $\infty \cdot
h(1,0)=h(c)$. Similarly,  $\infty \cdot h(1,1)=h(d)$. Moreover,
$h(1,0) + h(1,0) \neq h(1,0)$ and $h(1,1) + h(1,0) \neq h(1,1)$, so
$h(1,1)$ and $h(1,0)$ have empty infinite support and must be
incomparable elements. This contradicts the fact that $\infty \cdot
h(1,1)+\infty \cdot h(1,0)=\infty \cdot h(1,1)$. Therefore, $M$
cannot be isomorphic to a monoid given by equations.
\end{Proof}

Finally, we draw some consequences  for monoids of projective
modules of the results obtained in this section.

\begin{Cor} \label{fgwr} Let $R$ be a semilocal ring, let $\varphi \colon R\to
S$ be an onto ring homomorphism such that $\mathrm{Ker}\, \varphi
=J(R)$ and  $S\cong M_{n_1}(D_1)\times \cdots \times M_{n_k}(D_k)$
for suitable division rings $D_1,\dots ,D_k$. Assume that $\dim
_\varphi V^*(R)$ can be  defined by a system of inequalities such
that $\dim _\varphi V^*({}_RR)= D(\dim _\varphi V^*(R))$.

Then the monoids   $W(R)$, $W({}_RR)$, $V^*(R)$ and $V^*({}_RR)$ are
finitely generated. In addition, $W(R)$ and $W({}_RR)$ are
cancellative and unperforated.

If $P$ is a projective right module such that $\langle P\rangle \in
W(R)$ then $V(\mathrm{End}_R(P))$ is a cancellative, finitely
generated and unperforated monoid.
\end{Cor}

\begin{Proof} By Corollary \ref{wrvr} and Remark
\ref{finiteandinequalities}. The elements of $W(R)$ are the
solutions in $\No^k$ of the system of inequalities defining $M$. By
Proposition \ref{characterizations}, $W(R)$ is finitely generated
and unperforated. Being isomorphic to a submonoid of $\No^k$, $W(R)$
is also cancellative. The statement on $W({}_RR)$ follows by
symmetry.

By Lemma \ref{monoequality}, it follows that $V^*(R)$ and
$V^*({}_RR)$ are finitely generated.

Let $P$ be a projective right $R$-module such that $P/PJ(R)$ is
finitely generated. Using that the category of modules that are
direct summands of $P^n$, for some $n$,  is equivalent to the
category of finitely generated projective right modules over
$\mathrm{End}_R(P)$, we deduce that
\[V(\mathrm{End}_R(P))\cong \{x\in W(R)\mid \mbox{ there exists $n$
such that }x\le n\langle P\rangle\}=M\] Since $W(R)$ is finitely
generated, cancellative and unperforated then so is $M$.
\end{Proof}

\begin{Rem}
Observe that if $R/J(R)$ is right noetherian then $\langle P \rangle \in W(R)$
if and only if $P/PJ(R)$ is finitely generated. In this case $W(R)$ is finitely generated
whenever $V^{*}(R)$ is finitely generated.
\end{Rem}

For a general  semilocal ring we do not know whether the
endomorphism ring of a projective right $R$-module $P$ such that it
is finitely generated modulo the Jacobson radical must be again a
semilocal ring. We do not even know whether this happens for the
rings appearing in Theorem \ref{pbinequalities}. On the positive
side, Corollary \ref{fgwr} shows that, at least, the monoid
$V(\mathrm{End}_R(P))$ is of the \emph{correct type}, cf.
Proposition \ref{characterizations}.

\section{Realizing monoids defined by inequalities}\label{realizing}

We use the following result to construct semilocal rings with
prescribed $V^*(R)$.

\begin{Th} \label{milnorsemilocal} \emph{\cite{HP}} Let $R_1$
and $R_2$ be  semilocal rings, and let  $S= M_{m_1}(D'_1)\times
\cdots \times M_{m_\ell}(D'_\ell)$ for suitable division rings
$D'_1,\dots ,D'_\ell$. For $i=1,2$, let  $j_i\colon R_i \to S$ be
ring homomorphisms. Let $R$ be the ring that fits into the pullback
diagram
$$\begin{CD}
R_1@> j_1 >>S\\
@A i_1 AA @AA j_2 A\\
R@>> i_2> R_2
\end{CD}$$
Assume that $j_1$ is an onto ring homomorphism with kernel $J(R_1)$,
and that $J(R_2)\subseteq \mathrm{Ker} \, j_2$.  If $R_2/J(R_2)\cong
M_{n_1}(D_1)\times \cdots \times M_{n_k}(D_k)$ for $D_1,\dots ,D_k$
division rings, and $\pi \colon R_2\to M_{n_1}(D_1)\times \cdots
\times M_{n_k}(D_k)$ is an onto morphism with kernel $J(R_2)$ then
\begin{itemize}
\item[(i)] $i_2$ induces an onto ring homomorphism
$\overline{i_2}\colon R\to M_{n_1}(D_1)\times \cdots \times
M_{n_k}(D_k)$ with kernel $J(R)$. In particular, $R$ is a semilocal
ring and  $R/J(R)\cong R_2/J(R_2)$.
\item[(ii)] Let $\alpha \colon \dim _{\pi}  V^*(R_2)\to (\No^*)^{\ell}$ be the monoid homomorphism
induced by $j_2$. Then $$\dim _{\overline{i_2}} V^*(R)=\{x\in \dim
_{\pi }V^*(R_2)\mid \alpha (x)\in \dim _{ j_1}V^*(R_1)\},$$ and
$$\dim _{\overline{i_2}} V^*({}_RR)=\{x\in \dim
_{\pi }V^*({}_{R_2}R_2)\mid \alpha (x)\in \dim _{
j_1}V^*({}_{R_1}R_1)\}.$$
\end{itemize}
\end{Th}

\begin{Ex}\label{inequalities}  Let $k\in \N$, and let
$a_1,\dots,a_k, b_1,\dots ,b_k\in \No$. Let  $(n_1,\dots ,n_k)\in \N
^k$ be such that $a_1n_1+\cdots +a_kn_k= b_1n_1+\cdots +b_kn_k\in
\N$. For any field extension $F\subseteq F_1$, there exist a
 semilocal $F$-algebra $R$ and an onto morphism of
$F$-algebras $\varphi \colon R\to M_{n_1}(F_1)\times \cdots \times
M_{n_k}(F_1)$ with kernel $J(R)$ such that $\dim _{\varphi}V^*(R_R)$
is the set of solutions in $(\No^*)^k$ of the inequality
$a_1t_1+\cdots +a_kt_k\ge  b_1t_1+\cdots +b_kt_k$ and $\dim
_{\varphi}V^*({}_RR)$ is the set of solutions in $(\No^*)^k$ of the
inequality $a_1t_1+\cdots +a_kt_k\le  b_1t_1+\cdots +b_kt_k$.

Note that $\dim _\varphi  (\langle R\rangle )=(n_1,\dots ,n_k).$
\end{Ex}

\begin{Proof}Set $m=a_1n_1+\cdots +a_kn_k= b_1n_1+\cdots +b_kn_k$.

Let $T$ be a  semilocal $F$-algebra  with an onto algebra morphism
$j_1\colon T\to F_1\times F_1$ with $\mathrm{Ker}(j_1)=J(T)$, and
such that  $\mathrm{dim}_{j_1}\,V^*(T_T)=\{(x,y)\in (\No^*)^2\mid
x\ge y\}$ and $\mathrm{dim}_{j_1}\,V^* ({}_TT)=\{(x,y)\in
(\No^*)^2\mid y\ge x\}$. Such $T$ exists by Theorem~\ref{GS}. Let
$M_m(j_1)\colon M_m(T)\to M_m(F_1)\times M_m(F_1)$ be the induced
morphism.

Set $R_2=M_{n_1}(F_1)\times \dots \times M_{n_k}(F_1)$. Consider the
 morphism of $F$-algebras $j_2\colon R_2\longrightarrow M_{m}(F_1)\times
M_m(F_1)$ defined by
\[j_2(r_1,\dots ,r_k)=
\left( \left(\begin{array}{ccc}\smallmatrix r_1&\cdots &0\\
\vdots &\ddots\! ^{a_1)} &\vdots \\ 0&\cdots& r_1\endsmallmatrix
&\cdots &
0\\
&\ddots & \\
0&\cdots &\smallmatrix r_k&\cdots &0\\
\vdots&\ddots\! ^{a_k)} &\vdots \\ 0&\cdots& r_k\endsmallmatrix
\end{array}\right),\left(\begin{array}{ccc}\smallmatrix r_1&\cdots &0\\
\vdots &\ddots\! ^{b_1)} &\vdots \\ 0&\cdots& r_1\endsmallmatrix
&\cdots &
0\\
&\ddots & \\
0&\cdots &\smallmatrix r_k&\cdots &0\\
\vdots&\ddots\! ^{b_k)} &\vdots \\ 0&\cdots& r_k\endsmallmatrix
\end{array}\right)\right)
\]
The morphism $j_2$ induces the morphism of monoids $f\colon
(\No^*)^k\to \No^*\times \No^*$ defined by $f(x_1,\dots
,x_k)=(a_1x_1+\cdots +a_kx_k, b_1x_1+\cdots +b_kx_k)$. Hence,
$f(n_1,\dots ,n_k)=(m,m)$.

Let $R$ be the ring defined by the pullback diagram
$$\begin{CD}
M_m(T)@> M_m(j_1) >>M_{m}(F_1)\times M_m(F_1)\\
@A i_1 AA @AA j_2 A\\
R@>> \varphi > M_{n_1}(F_1)\times \cdots \times M_{n_k}(F_1)
\end{CD}$$

 Applying
Theorem~\ref{milnorsemilocal} (i), we  conclude  that $R$ is a
 semilocal $F$-algebra
 and that $\varphi$ is an onto morphism of $F$-algebras with kernel $J(R)$.
By Theorem~\ref{milnorsemilocal}(ii),
 $(x_1,\dots ,x_k)\in \dim _{\varphi}V^*(R_R)$ if and only
if $f(x_1,\dots ,x_k)\in \dim _{M_m(j_1)}V^*(M_m(T))$ if and only if
$a_1x_1+\dots +a_kx_k\ge b_1x_1+\cdots+b_kx_k$. Similarly,
$(x_1,\dots ,x_k)\in \dim _{\varphi}V^*({}_RR)$ if and only if
$a_1x_1+\dots +a_kx_k\le b_1x_1+\cdots+b_kx_k$.
\end{Proof}

Now we are ready to prove Theorem \ref{pbinequalities}.

\begin{Proofpbinequalities} Let  $M$ be the monoid defined by the
system of inequalities,
\[D\left(\begin{array}{c}t_1\\\vdots \\ t_k\end{array}\right)\in \left(\begin{array}{c}m_1\No^*\\
\vdots \\ m_n\No^*\end{array} \right)\quad(*)\qquad \mbox{and}
\qquad
E_1\left(\begin{array}{c}t_1\\\vdots \\
t_k\end{array}\right)\le E_2\left(\begin{array}{c}t_1\\\vdots \\
t_k\end{array}\right)\quad(**)\] where  $D\in M_{n\times k}(\No)$,
$E_1,$ $E_2\in M_{\ell \times k}(\No)$, $n,\ell\ge 0$ and $m_1,\dots
,m_n\in \N$ , $m_i\ge 2$ for any $i\in \{1,\dots ,n\}$.

By \cite[Theorem 5.3]{HP} we know the following.

\noindent\textbf{Step 1.} \emph{There exist  a field extension $E$
of $F$, a (noetherian) semilocal $F$-algebra $R_1$ and an onto
morphism of $F$-algebras $\varphi _1\colon R_1\to M_{n_1}(E)\times
\cdots \times M_{n_k}(E)$ such that $\dim _{\varphi _1} V^*(R_1)$ is
the set of solutions in $(\No^*)^k$ of the system of congruences
$(*)$.}

Now we need to prove,

\noindent\textbf{Step 2.} \emph{There exist a   semilocal
$F$-algebra $R_2$ and an onto morphism of $F$-algebras $\varphi
_2\colon R_2\to M_{n_1}(E)\times \cdots \times M_{n_k}(E)$ such that
$\dim _{\varphi _2} V^*(R_2)$ is the set of solutions in $(\No^*)^k$
of the system of inequalities $(**)$ and $\dim _{\varphi _2}
V^*(_{R_2}R_2)$ is the set of solutions in $(\No^*)^k$ of the system
of inequalities $D(**)$.}

If $\ell=0$, that is, if $(**)$ is empty we set
$R_2=M_{n_1}(E)\times \cdots \times M_{n_k}(E)$ and $\varphi
_2=\mathrm{Id}$. Assume $\ell >0$. Therefore,  we can assume that
none of the rows in $E_1$ and, hence, in $E_2$ are zero.

By Example~\ref{inequalities}, for $i=1,\dots ,\ell$, there exist a
noetherian semilocal $F$-algebra $T_i$ and an onto morphism of
$F$-algebras $\pi _i\colon T_i \to M_{n_1}(E)\times \dots \times
M_{n_k}(E)$ with kernel $J(T_i)$ and such that $\dim _{\pi
_i}V^*(T_i)$ is the set of solutions in $(\No^*)^k$ of the $i$-th
inequality defined by the matrices $E_1$ and $E_2$, and  $\dim _{\pi
_i}V^*({}_{T_i}T_i)$ is the set of solutions in $(\No^*)^k$ of the
reversed inequality.

Let $R_2$ be the pullback of $\pi _i$, $i=1,\dots ,\ell$. By
Theorem~\ref{milnorsemilocal}, $R_2$ is a  semilocal $F$-algebra
with an onto morphism of $F$-algebras $\varphi_2\colon R_2\to
M_{n_1}(E)\times \dots \times M_{n_k}(E)$ with kernel $J(R_2)$.
Moreover, $\dim _{\varphi _2} V^*(R_2)$ is the set of  solutions of
the inequalities $(**)$ and $\dim _{\varphi _2} V^*(_{R_2}R_2)$ is
the set of  solutions of the inequalities $D(**)$. This concludes
the proof of Step 2.

\medskip

Finally, set $R$ to be the pullback of $\varphi _i \colon R_i\to
M_{n_1}(E)\times \dots \times M_{n_k}(E)$, $i=1,2$. By
Theorem~\ref{milnorsemilocal}, $R$ is a semilocal $F$-algebra with
an onto  morphism of $F$-algebras $\varphi\colon R\to
M_{n_1}(E)\times \dots \times M_{n_k}(E)$ with kernel $J(R)$. The
elements in $\dim _{\varphi } V^*(R_R)$ are the solutions of $(*)$
and $(**)$,  and the ones in $\dim _{\varphi } V^*({}_RR)$ are the
elements of $D(M)$.

The description of the images via $\dim _\varphi$ of $V(R)$, $W(R)$
and $W({}_RR)$ follows from Remark \ref{finiteandinequalities}.
\end{Proofpbinequalities}

\renewcommand{\baselinestretch}{1}

\end{document}